\documentclass[12pt,fullpage]{article}

\usepackage{amsfonts,graphicx,amsmath,xr,epsfig,verbatim,colortbl}

\def\N{\mathbb{N}}
\def\R{\mathbb{R}}
\def\T{\mathbb{T}}

\def\epsilon{\varepsilon}
\def\tilde{\widetilde}
\newcommand{\Bl}{\color{blue}}

\def\eps{\epsilon}

\newcommand{\be}{\begin{equation}}
\newcommand{\ee}{\end{equation}}
\newcommand{\baa}{\begin{array}}
\newcommand{\eaa}{\end{array}}
\newcommand{\ba}{\begin{eqnarray}}
\newcommand{\ea}{\end{eqnarray}}

\newtheorem{theo}{\bf Theorem}[section]
\newtheorem{lem}[theo]{\bf Lemma}
\newtheorem{pro}[theo]{\bf Proposition}

\newtheorem{rem}[theo]{\bf Remark}

\begin{document}
\date{}
\title{\bf{Speed-up of combustion fronts in shear flows}}
\author{Fran\c cois Hamel$^{\hbox{\small{ a}}}$$\ $ and Andrej Zlato{\v{s}}$^{\hbox{\small{ b }}}$
\\
\footnotesize{$^{\hbox{a }}$Aix-Marseille Universit\'e \& Institut Universitaire de France}\\
\footnotesize{LATP, Facult\'e des Sciences et Techniques, F-13397 Marseille Cedex 20, France}\\
\footnotesize{$^{\hbox{b }}$Department of Mathematics, University of Wisconsin, Madison, WI 53706, USA}
}
\maketitle

\begin{abstract} 
This paper is concerned with the analysis of speed-up of reaction-diffusion-advection traveling fronts in infinite cylinders with periodic boundary conditions. The advection is a shear flow with a large amplitude and the reaction is nonnegative, with either positive or zero ignition temperature. The unique or minimal speeds of the traveling fronts are proved to be asymptotically linear in the flow amplitude as the latter goes to infinity, solving an open problem from \cite{b}. The asymptotic growth rate is characterized explicitly as the unique or minimal speed of traveling fronts for a limiting degenerate problem, and the convergence of the regular traveling fronts to the degenerate ones is proved for positive ignition temperatures under an additional H\" ormander-type condition on the flow.
\end{abstract}

%%%%%%%%%%%%%%%%%%%%%%%%%%%%%%%%%%%%%%%%
%%%%%%%%%%%%%%%%%%%%%%%%%%%%%%%%%%%%%%%%

\section{Introduction and main results}\label{intro}

In this paper we analyze the asymptotic limit of traveling fronts and of their speeds for reaction-diffusion equations in the presence of strong shear flows. More precisely, we consider the model
\be\label{eqv}
v_t+A\,\alpha(y)\,v_x=\Delta v+f(v),\quad t\in\R,\ (x,y)\in\R\times\R^{N-1}
\ee
when $A\gg 1$. The physical dimension~$N\ge 1$ is arbitrary. The {\it reaction function} $f:[0,1]\to\R$ is assumed to be of class $C^{1,\delta}([0,1])$ for some $\delta>0$ and to satisfy 
\be\label{f}
\exists\,\theta\in[0,1):\quad f=0\hbox{ on }[0,\theta]\cup\{1\},\quad f>0\hbox{ on }(\theta,1), \quad f \text{ is non-increasing near 1}.
\ee
The function~$\alpha:\R^{N-1}\to\R$ is assumed to be of class~$C^{1,\delta}(\R^{N-1})$ and $(1,\ldots,1)$-periodic  (for the sake of simplicity; other periods are handled identically) so that $\alpha\in C^{1,\delta}(\T^{N-1})$. The coefficient $A\,\alpha(y)$ is then the $x$-component of an incompressible {\it shear flow}~$A\,q(x,y)$ with {\it amplitude}~$A$ and {\it flow profile}
$$q(x,y)=(\alpha(y),0,\ldots,0).$$
We are interested in the asymptotic strong-flow regime $A\to+\infty$.

Equation \eqref{eqv} arises in models of flame propagation, especially when $0<\theta<1$, and the quantity~$0\le v\le 1$ then stands for  normalized temperature (see~\cite{k2}). The real number~$\theta$ is then the {\it ignition temperature}, below which the reaction is suppressed. When $\theta=0$, this equation also arises in chemical, biological, and ecological models, and~$v$ typically stands for the density of a substance or a species (see~\cite{f,mu,sk}).

We will be interested in {\it traveling fronts} for \eqref{eqv}, that is, solutions of the form
$$v(t,x,y)=u(x-ct,y).$$
Here the {\it front profile} $u$ is $(1,\ldots,1)$-periodic in $y\in\R^{N-1}$ and $c\in\R$ is the {\it front speed}. These solutions, which are invariant in the frame moving with speed~$c$ along the $x$-axis, play an essential role in the study of large-time behavior of the processes modeled by~(\ref{eqv}). 

If $0<\theta<1$, then it is known \cite{bh1,Xin}  that for each amplitude~$A\in\R$ there exists a {\it unique speed}~$c=c^*(Aq,f)$ and a unique (up to shifts in $x$) profile~$u:\R\times\T^{N-1}\to(0,1)$, which connects~$0$ to $1$. That is~$u$ is a classical solution of the elliptic equation
\be\label{equ}\left\{\baa{rcl}
\Delta u+(c-A\,\alpha(y))\,u_x+f(u) & \!\!=\!\! & 0\ \hbox{ in }\R\times\T^{N-1},\vspace{3pt}\\
0 \ \equiv \ u(+\infty,\cdot)\  <\ u\ <\ u(-\infty,\cdot) & \!\!\equiv\!\! & 1\ \hbox{ in }\R\times\T^{N-1},\eaa\right.
\ee
where the above limits $0$ and $1$ are uniform in $\T^{N-1}$. Moreover,  $u$ is decreasing in~$x$. 

On the other hand, if $\theta=0$, then there is a {\it minimal speed} $c^*(Aq,f)$ such that problem~(\ref{equ}) has a solution~$u$ with speed $c$ if and only if $c\ge c^*(Aq,f)$. Moreover, for each such $c$ there is a decreasing-in-$x$ solution, and all solutions are such if $f'(0)>0$ \cite{bh1}. 

The question of existence of traveling fronts for shear flows was first considered in~\cite{bll,bn2} in infinite cylinders with bounded cross sections and Neumann boundary conditions, and the case with periodic boundary conditions in~$y$ can be treated similarly.  Since then, a considerable amount of research has been directed at the problem, in particular, first results on stability of fronts and long-time convergence of large classes of solutions to them appeared in~\cite{aw,blr,k2,r1,r2}.  Existence of fronts in the case of general, not necessarily shear, flows which are periodic in all spatial variables (on periodic domains in $\R^N$) has also been studied. These {\it pulsating  fronts} are of the type $v(t,z)=U(z\cdot e-ct,z)$, with~$z=(x,y)$, $e\in\R^N$ a unit vector, and~$U$ periodic in the last $N$ arguments $z=(z_1,\ldots,z_N)\in\R^N$, and their existence was first proved in~\cite{bh1,x1}. We also refer to~\cite{e2,h1,hps} for variational min-max type formulas for the unique and minimal propagation speeds of traveling and pulsating fronts, and to~\cite{bhbook,x2} for many further references on propagation phenomena for reaction-diffusion equations.

The main question addressed in this paper is related to the behavior of the traveling fronts~$u$ of~(\ref{equ}) and the propagation speeds~$c^*(Aq,f)$ in strong shear flows, that is as the amplitude~$A$ goes to $+\infty$. We especially aim at quantifying the speed-up induced by the underlying strong advection. Roughly speaking, in combustion models, a stronger advection is going to increase the width of the reaction zone, and hence enhance the propagation speed. It is clear that if $\alpha\le\beta$, then, for each $A\ge 0$, the corresponding unique or minimal speeds associated with~$\alpha$ and~$\beta$ are ordered, that is $c^*(Aq_{\alpha},f)\le c^*(Aq_{\beta},f)$, where $q_{\alpha}(x,y)=(\alpha(y),0,\ldots,0)$ and $q_{\beta}(x,y)=(\beta(y),0,\ldots,0)$, see~\cite{bhbook,bn2}. However, in general the speeds for flows $A\,\alpha(y)$ and $A'\,\alpha(y)$ cannot be compared if~$A \neq A'$. Therefore, the average effect of the amplitude~$A$ and the limiting behavior as $A\to+\infty$ are a priori not clear. Nevertheless, the comparison principle immediately shows
\be\label{limsup}
\limsup_{A\to+\infty}\frac{c^*(Aq,f)}{A}\le\max_{\T^{N-1}}\,\alpha
\ee
(see~\cite{bhbook,bn2}). Furthermore, it also follows from~\cite{kr} that
\be\label{liminf}
\liminf_{A\to+\infty}\frac{c^*(Aq,f)}{A}\ge \int_{\T^{N-1}} \alpha(y)dy,
\ee
with equality only if $\alpha$ is constant.  (In \cite{kr} only mean-zero $\alpha$ were considered, with the above bound being 0, but the general case follows from $c^*(Aq,f) = c^*(A(q-p),f)+A\beta$ if $p(x,y)=(\beta,0,\cdots,0)$.) That is, the front speed in the presence of a strong shear flow is between two linear functions of the amplitude of the flow.  The main result of this paper shows that this speed-up is truly {\it asymptotically linear} as $A\to+\infty$, that is,  $\lim_{A\to +\infty} c^*(Aq,f)/A$ exists.   Furthermore, we characterize this limit, in both cases $0<\theta<1$ and $\theta=0$, in terms of solutions to a second-order degenerate equation. 

We start with the positive ignition temperature  case  $0<\theta<1$.  The existence of the limit $\lim_{A\to +\infty} c^*(Aq,f)/A$ in this case was stated as an open problem in \cite{b}, where it was proved for KPP reactions (see below).

\begin{theo}\label{th1}
Assume $(\ref{f})$ with $0<\theta<1$. Then there exists~$\gamma^*(q,f)\ge \int_{\T^{N-1}} \alpha(y)dy$ (with equality only if $\alpha$ is constant) such that the unique speeds $c^*(Aq,f)$ of the problem~$(\ref{equ})$ satisfy
\be\label{gamma}
\lim_{A\to +\infty} \frac{c^*(Aq,f)}{A} = \gamma^*(q,f).
\ee
Furthermore, there exists $U\in L^{\infty}(\R\times\T^{N-1})$ such that $\nabla_yU\in L^2(\R\times\T^{N-1})\cap L^{\infty}(\R\times\T^{N-1})$ and $U$ is a (distributional) solution of
\be\label{eqU}\left\{\baa{ll}
\Delta_yU+(\gamma-\alpha(y))\,U_x+f(U)=0 & \hbox{in }\mathcal{D}'(\R\times\T^{N-1}),\vspace{3pt}\\
0\le U\le 1 & \hbox{a.e. in }\R\times\T^{N-1},\vspace{3pt}\\
\lim_{x\to+\infty} U(x,y)\equiv 0 & \text{uniformly in }\T^{N-1},\vspace{3pt}  \\
\lim_{x\to-\infty} U(x,y)\equiv 1 & \text{uniformly in }\T^{N-1},
\eaa\right.
\ee
with $\gamma =  \gamma^*(q,f)$, and each sequence $A_n\to\infty$ has a subsequence along which the functions $U_{A_n}(x,y)=u_{A_n}(A_n x,y)$ (with $u_{A_n}$ the unique solution of \eqref{equ} for $A=A_n$, translated in $x$ so that  $\max_{y\in \T^{N-1}}u_{A_n}(0,y)=\theta$) converge a.e. to a solution of \eqref{eqU} with $\gamma =  \gamma^*(q,f)$ and $U$ non-increasing in $x$.
Finally, if  there is $r>0$ such that
\be\label{hypalpha}
\alpha\in C^{\infty}(\T^{N-1}) \qquad \text{and} \qquad
\sum_{1\le |\zeta|\le r} |D^\zeta \alpha(y)|>0 \text{ for all $y\in\T^{N-1}$,} \footnote{For $\zeta=(\zeta_1,\ldots,\zeta_{N-1})\in\N^{N-1}$, we let $|\zeta|=\zeta_1+\cdots+\zeta_{N-1}$ and $D^{\zeta}\alpha(y)=\frac{\partial^{|\zeta|}\alpha}{\partial y_1^{\zeta_1}\cdots\partial y_{N-1}^{\zeta_{N-1 }}}(y)$.}
\ee
then  the pair~$(\gamma,U)$ solving~$(\ref{eqU})$ is unique up to shifts of $U$ in~$x$, with $U\in C^{1,\delta}(\R\times\T^{N-1})$ satisfying  $0<U<1$ and $ U_x<0$ on $\R\times\T^{N-1}$.
\end{theo}

\begin{rem} {\rm
This result continues to hold when the ignition reaction $f$  is only Lipschitz and $\alpha\in C^{\delta}(\T^{N-1})$ for some $\delta>0$.  Indeed, our proof does not use the extra smoothness, and \cite[Theorem 1.6]{ZlaGenfronts} shows that there still exists a unique solution $u$ to \eqref{equ}, which is decreasing in $x$ due to being  the limit of solutions with the same property and due to the strong maximum principle applied to $u_x$.
}\end{rem}

\begin{rem} {\rm
The condition \eqref{hypalpha} ensures that the differential operator $L=\Delta_y + (\gamma-\alpha(y)) \partial_x$ from \eqref{eqU} satisfies H\" ormander's hypoellipticity condition for any $\gamma\in\R$, namely that  the smooth vector fields $\partial_{y_1},\dots,\partial_{y_{N-1}}, (\gamma-\alpha(y)) \partial_x$ and their commutators up to order $r+1$ span all of $\R^N$ at any point $(x,y)$.  Because of this, we will be able to apply the results from \cite{rs} to questions of regularity of solutions of \eqref{eqU}.
}\end{rem}

\begin{rem} {\rm
It remains an open question whether  uniqueness of the pair $(\gamma,U)$ solving \eqref{eqU} as well as $0<U<1$ and $ U_x<0$ hold for all non-constant $\alpha$. Note that any $y$-independent not necessarily continuous $U$, taking values in $[0,\theta]\cup\{1\}$ and having the prescribed limits, along with $\gamma=\alpha$, is a solution for constant $\alpha$.
}\end{rem}

On the other hand, for zero ignition temperature $\theta=0$ the following conclusion holds.  We define for ${\theta'}\in(0,1/4]$ the ``cut-off'' ignition reaction $f_{\theta'}(u):= f(u)\chi(u/{\theta'})$, where $\chi:\R\to [0,1]$ is a smooth non-decreasing function with $\chi(v)=0$ if $v\le 1$ and $\chi(v)=1$ if $v\ge 2$.

\begin{theo}\label{th2}
Assume~$(\ref{f})$ with $\theta=0$ and also $f\in C^{1,1}([0,1])$. 
Then the limit
$$\gamma^*(q,f):=\lim_{{\theta'}\to 0} \gamma^*(q,f_{{\theta'}})$$
exists (with $\gamma^*(q,f_{{\theta'}})$ from Theorem \ref{th1}) and~$(\ref{gamma})$ holds for the minimal speeds $c^*(Aq,f)$ of pro\-blem~$(\ref{equ})$. Furthermore, for any $\gamma\ge\gamma^*(q,f)$, there exists a solution~$U$ of~$(\ref{eqU})$ such that  $\nabla_yU\in L^2(\R\times\T^{N-1})\cap L^{\infty}(\R\times\T^{N-1})$ and $U$ is non-increasing in $x$. Lastly, if~$(\ref{hypalpha})$ holds, then all solutions~$U$ of~$(\ref{eqU})$ are classical and satisfy $\gamma\ge\gamma^*(q,f)$.
\end{theo}

When $\theta=0$ in~(\ref{f}) and the nonlinearity $f$ is of the Kolmogorov-Petrovsky-Piskunov (KPP)~\cite{kpp} type, that is $f(u)\le f'(0)u$ for all $u\in[0,1]$, then the limit $\gamma^*(q,f)$ given in Theorem~\ref{th2} can be expressed as
\be\label{gamma*kpp}
\gamma^*(q,f)=\max_{\substack{w\in H^1(\T^{N-1})\backslash\{0\}\\ \|\nabla w\|_{L^2(\T^{N-1})}^2\le f'(0)\|w\|_{L^2(\T^{N-1})}^2}}\frac{\displaystyle{\int_{\T^{N-1}}}\alpha(y)\,w(y)^2 dy}{\displaystyle{\int_{\T^{N-1}}}w(y)^2 dy}
\ee
(see~\cite{b,he,z5}). Using this, one can conclude \cite{he,z5}
$$\displaystyle{\mathop{\lim}_{M\to+\infty}}\left(\displaystyle{\mathop{\lim}_{A\to+\infty}}\ \displaystyle{\frac{c^*(Aq,Mf)}{A}}\right)=\max_{y\in \T^{N-1}}\alpha(y)$$
and
$$\displaystyle{\mathop{\lim}_{M\to0^+}}  \frac 1{\sqrt{M}} \left(\displaystyle{\mathop{\lim}_{A\to+\infty}}\!\displaystyle{\frac{c^*(Aq,Mf)}{A} -\int_{\T^{N-1}} \alpha(y) dy}\right)\!=\! 2\,\sqrt{f'(0)}\sup_{\substack{w\in H^1(\T^{N-1})\backslash\{0\} \\ \int_{\T^{N-1}} w(y) dy =0}}\!\!\frac{\displaystyle{\int_{\T^{N-1}}}\!\!\alpha(y)w(y)dy}{\|\nabla w\|_{L^2(\T^{N-1})}}.$$
The first of these limits does not depend on $f$ and so it holds for any $f$ with $f'(0)>0$.  Indeed, then we can find two KPP reactions $g,h$ such that $g\le f \le h$, and the claim follows from the fact that the limits for $g,h$ are equal and from the well known relation $c^*(Aq,g)\le c^*(Aq,f) \le c^*(Aq,h)$~\cite{bn2}.

Similar results hold for more general heterogeneous equations with periodic coefficients and KPP reactions. We refer to~\cite{bmr,bhn2,bhr2,h2,hr,mr,we} for further results on existence and qualitative properties of KPP traveling and pulsating fronts and to~\cite{abp,b,bhn1,e1,ek,rz,z5} for  results on the asymptotics of the minimal speeds in strong shear and periodic flows in the KPP case.  

We also note that for cellular flows in two spatial dimensions, the speed-up is known to be  $O(A^{1/4})$ rather than linear (see~\cite{abp,ckor,kr,nr,rz,z4}), and characterization of those two-dimensional periodic $q$ for which $c^*(Aq,f)$ is unbounded in $A$ (i.e., speed-up occurs) was  obtained in~\cite{rz,z4}.  In fact, these two-dimensional results hold for general reactions \cite{z4}.

\begin{rem}{\rm In the case of KPP reactions and shear flows with $\int_{\T^{N-1}} \alpha(y)dy\ge 0$, it is also known that $c^*(Aq,f)$ is non-decreasing with respect to $A\ge 0$ and $c^*(Aq,f)/A$ is non-increasing with respect to $A\ge 0$, see~\cite{b,bhbook}. However, the monotonicity of the quantities $c^*(Aq,f)$ and $c^*(Aq,f)/A$ for general non-KPP reactions~$f$ when either $\theta=0$ or $0<\theta<1$, is still not known. These questions are much more intricate than in the KPP case, due to the lack of explicit formulas for the propagation speeds. However,  Theorems~\ref{th1} and~\ref{th2} show that the limit $\gamma^*(q,f)$ of the normalized speeds $c^*(Aq,f)/A$ can still be identified as the unique or minimal speed of fronts for a degenerate elliptic problem.}
\end{rem}

We thank  Tom Kurtz and Daniel Stroock for useful discussions and pointers to references.  FH is indebted to the Alexander von~Humboldt Foundation for its support.  His work was also supported by the French {\it Agence Nationale de la Recherche} through the project PREFERED.  AZ was supported in part by NSF grants DMS-1113017 and DMS-1056327, and by an Alfred P. Sloan Research Fellowship.  Part of this work was carried out during visits by FH to the Departments of Mathematics of the Universities of Chicago and Wisconsin and by AZ to the Facult\'e des Sciences et Techniques, Aix-Marseille Universit\'e, the hospitality of which is gratefuly acknowledged. 

%%%%%%%%%%%%%%%%%%%%%%%%%%%%%%%%%%%%%%%%
%%%%%%%%%%%%%%%%%%%%%%%%%%%%%%%%%%%%%%%%

\section{Existence of solutions $(\gamma,U)$ of~(\ref{eqU}) for $\theta\in[0,1)$}\label{sec2}

In this section, we prove the existence of solutions $(\gamma,U)$ of~(\ref{eqU}) as claimed in Theorems \ref{th1} and \ref{th2}. Such solutions are obtained as the $A\to+\infty$ limit of the solutions $(c_A,u_A)$ of~(\ref{equ}), after a suitable normalization and a scaling in $x$. The scaling, given in~(\ref{defUA}) below, makes the first-order coefficient bounded but it makes the equation degenerate as $A\to+\infty$. The degeneracy of the limiting equation makes the analysis of the properties of the limiting solutions more complicated than in the usual regular case.

We will assume in this section that
\be\label{2.1}
\int_{\T^{N-1}}\alpha(y)dy=0,
\ee
which we can do without loss of generality because if $p(x,y)=(\beta,0,\cdots,0)$, then $c^*(Aq,f) = c^*(A(q-p),f)+A\beta$ and the two $u_A$ are identical.  Moreover, we will also assume that $\alpha$ is not identically $0$. Otherwise, $c^*(Aq,f)=c^*(0,f)$ does not depend on $A$ and the existence results in both theorems are immediate with $\gamma^*(q,f)=0$ and  $U(x,y)=\chi_{(-\infty,0)}(x)$ (for $\gamma=0$) or $U(x,y)=V(x)$ such that $V_x=-f(V)/\gamma$, $V({\Bl +}\infty)=0$, and $V(-\infty)=1$ (for $\gamma>0$ in Theorem \ref{th2}), which obviously exists. 

Set
$$\gamma^*_A=\frac{c^*(Aq,f)}{A} >0$$
for each $A\ge 1$, where $c^*(Aq,f)$ is the unique (resp.~minimal) speed of traveling fronts for~(\ref{equ}) when $0<\theta<1$ (resp.~$\theta=0$). Notice first that, because of~(\ref{limsup}) and~(\ref{liminf}), and since the speeds~$c^*(Aq,f)$ are positive and continuous in~$A$ (see~\cite{bhbook,bn2}), it follows that
\be\label{infsup}
 0<\underline{\gamma}:=\inf_{A\ge 1}\,\gamma^*_A\le\sup_{A\ge 1}\,\gamma^*_A =:\overline{\gamma}<+\infty.
\ee\par
Now, let $(c_A,u_A)_{A\ge 1}$ be a family of solutions of~(\ref{equ}), and define
$$\gamma_A=\frac{c_A}{A}.$$
We will only consider solutions with $u_{A,x}< 0$ {(which is always possible from~\cite{bh1})} and $\gamma_A$ uniformly bounded in $A\ge 1$. Notice that by~(\ref{infsup}), the family~$(\gamma_A)_{A\ge 1}$ is automatically bounded from below by a positive constant, and  it is also bounded from above when $0<\theta<1$.   For each $A\ge 1$, we set
\be\label{defUA}
U_A(x,y)=u_A(Ax,y)\ \hbox{ for all }(x,y)\in\R\times\T^{N-1}.
\ee
The functions $U_A$ satisfy
\be\label{eqUA}\left\{\baa{rcl}
\Delta_yU_A+A^{-2}U_{A,xx}+(\gamma_A-\alpha(y))\,U_{A,x}+f(U_A) & \!\!=\!\! & 0\ \hbox{ in }\R\times\T^{N-1},\vspace{3pt}\\
0 \ <\ U_A\  & \!\! < \!\! & 1\ \hbox{ in }\R\times\T^{N-1},\vspace{3pt} \\
U_{A,x} & \!\!<\!\! & 0 \ \hbox{ in }\R\times\T^{N-1},\vspace{3pt}\\
 U_A(+\infty,\cdot)\  \equiv \ 0, \quad  U_A(-\infty,\cdot) & \!\!\equiv\!\! & 1\ \hbox{ uniformly in }\T^{N-1}.\eaa\right.
\ee
By standard elliptic estimates, each function $U_A$ is of class $C^{3,\delta}(\R\times\T^{N-1})$.

We first derive some uniform integral and pointwise gradient estimates for the functions~$U_A$.

\begin{lem}\label{lem1} For any family $(c_A,u_A)_{A\ge 1}$ of solutions of~$(\ref{equ})$ such that $u_{A,x}< 0$ and $(c_A/A)_{A\ge 1}$ is bounded, the functions $U_A$ defined by~$(\ref{defUA})$ satisfy
$$\sup_{A\ge1}\Big(\|U_{A,x}\|_{L^1(\R\times\T^{N-1})}+\|\nabla_yU_A\|_{L^2(\R\times\T^{N-1})}+\|\nabla_yU_A\|_{L^{\infty}(\R\times\T^{N-1})}\Big)<+\infty.$$
\end{lem}

\noindent{\bf{Proof.}} First, since $U_{A,x}<0$ in $\R\times\T^{N-1}$ and $U_A(+\infty,\cdot)=0$, $U_A(-\infty,\cdot)=1$, it follows immediately that~$\|U_{A,x}\|_{L^1(\R\times\T^{N-1})}=1$.  Next integrate~(\ref{eqUA}) over $(-M,M)\times\T^{N-1}$ and pass to the limit as $M\to+\infty$. Since $U_{A,x}(\pm\infty,\cdot)=0$ from standard elliptic estimates, since~$\alpha$ has zero average and since $f\ge 0$ in~$[0,1]$, we get that~$f(U_A)\in L^1(\R\times\T^{N-1})$ and
$$\int_{\R\times\T^{N-1}}f(U_A(x,y))\,dx\,dy=\gamma_A$$
for each $A\ge 1$. Now multiply~(\ref{eqUA}) by $U_A$, integrate over $(-M,M)\times\T^{N-1}$ and pass to the limit as $M\to+\infty$. It follows that
$$\int_{\R\times\T^{N-1}}\Big(|\nabla_yU_A|^2+A^{-2}U_{A,x}^2\Big)\,dx\,dy=\int_{\R\times\T^{N-1}}f(U_A(x,y))\,U_A(x,y)\,dx\,dy-\frac{\gamma_A}{2}\le\frac{\gamma_A}{2}.$$
Therefore, the assumption of boundedness of the family $(\gamma_A)_{A\ge 1}$ yields the boundedness of the family $(\|\nabla_yU_A\|_{L^2(\R\times\T^{N-1})})_{A\ge 1}$.

Lastly, since~$U_A\in C^3(\R\times\T^{N-1})$ and the fa\-mily~$(\gamma_A)_{A\ge 1}$ is assumed to be bounded, it then follows from~\cite{bh2} that
$$\sup_{A\ge 1}\,\|\nabla_yU_A\|_{L^{\infty}(\R\times\T^{N-1})}<+\infty,$$
which completes the proof of Lemma~\ref{lem1}.\hfill$\Box$\break

We next prove the claims in the second sentence of Theorem \ref{th1} when $0<\theta<1$ with some $\gamma>0$ which is the limit of a subsequence of $c^*(A_nq,f)/A_n$.  The claim as stated will be established when later we prove \eqref{gamma}.

\begin{pro}\label{proex1} Assume $0<\theta<1$. Then any $A_n\to+\infty$ has a subsequence along which the (unique) solution $U_{A_n}$  of \eqref{eqUA}, normalized as in Theorem \ref{th1}, converges a.e. to $U$ and $\gamma_{A_n}$ converges to $\gamma$, with $(\gamma,U)$ a solution of~$(\ref{eqU})$ such that $0<\gamma\le\max_{\T^{N-1}}\,\alpha$ and $\nabla_yU\in L^2(\R\times\T^{N-1})\cap L^{\infty}(\R\times\T^{N-1})$.
\end{pro}

\noindent{\bf{Proof.}} {\it Step 1: general properties of the limit function $U$.} The assumptions of Lemma \ref{lem1} are satisfied since $0<\theta<1$.
Because of~(\ref{limsup}) and~(\ref{infsup}), there is a subsequence of $A_n$ (which we again call $A_n$) and  a real number~$\gamma$ such that
\be\label{gamma+}
0<\gamma\le\max_{\T^{N-1}}\,\alpha\ \hbox{ and }\ \gamma_{A_n}\to\gamma\hbox{ as }n\to+\infty.
\ee
On the other hand,  $0<U_{A_n}<1$ in $\R\times\T^{N-1}$ and the sequence~$(U_{A_n})_{n\in\N}$ is  bounded in $W^{1,1}_{\rm loc}(\R\times\T^{N-1})$ by Lemma~\ref{lem1}. Up to extraction of another subsequence, there exists then a function $U\in L^{\infty}(\R\times\T^{N-1})$ such that $U_{A_n}\to U$ almost everywhere in~$\R\times\T^{N-1}$ as $n\to+\infty$. Furthermore, $0\le U\le 1$ and
\be\label{monotone}
U(\cdot+h,\cdot)\le U\hbox{ a.e. in }\R\times\T^{N-1}
\ee
 for all $h\ge 0$.  Similarly, $\nabla_yU\in L^2(\R\times\T^{N-1})\cap L^{\infty}(\R\times\T^{N-1})$ by Lemma~\ref{lem1}. The equation
\be\label{eqUbis}
\Delta_yU+(\gamma-\alpha(y))\,U_x+f(U)=0
\ee
also holds in the sense of distributions in $\R\times\T^{N-1}$, by passing to the limit as $n\to+\infty$ in~(\ref{eqUA}), with $A=A_n$.

{\it Step 2: the limit function $U$ has constant limits as $x\to\pm\infty$.} Take any non-decreasing sequence $(x_n)_{n\in\N}$ converging to $+\infty$. The sequences of functions $(U^{\pm}_n)_{n\in\N}$ defined by
$$U^{\pm}_n(x,y)=U(x\pm x_n,y)\ \hbox{ for }(x,y)\in\R\times\T^{N-1}$$
are monotone and bounded in $L^{\infty}(\R\times\T^{N-1})$. They converge almost everywhere in $\R\times\T^{N-1}$ and in $L^p((-M,M)\times\T^{N-1})$, for all $M>0$ and $1\le p<+\infty$, to two functions $U^{\pm}_{\infty}$ such that $0\le U^{\pm}_{\infty}\le 1$ a.e. in~$\R\times\T^{N-1}$. Furthermore,
$$U^{\pm}_{\infty}(x+h,y)=U^{\pm}_{\infty}(x,y)\hbox{ a.e. in }\R\times\T^{N-1}\hbox{ for all }h\in\R,$$
due to~(\ref{monotone}) and $x_n\to+\infty$, and $\nabla_yU^{\pm}_{\infty}=0$ since $\nabla_yU\in L^2(\R\times\T^{N-1})$. In other words, the functions $U^{\pm}_{\infty}$ are constant. Since the equation~(\ref{eqUbis}) is satisfied by each~$U^{\pm}_n$ and thus by the limiting (constant) functions~$U^{\pm}_{\infty}$ in the sense of distributions, one concludes that
\be\label{fuinfty}
f(U^{\pm}_{\infty})=0.
\ee
By monotonicity in $x$, the limits $U^{\pm}_{\infty}$ do not depend on the sequence~$(x_n)_{n\in\N}$. Therefore, the limits $U(\pm\infty,\cdot)$ exist and are two constants belonging to the set of zeros of $f$ in $[0,1]$.

{\it Step 3: the limit function $U$ has uniform limits $0$ and $1$ as $x\to\pm\infty$.} The last thing to do is to prove that $U(+\infty,\cdot)=0$ and $U(-\infty,\cdot)=1$ uniformly in $y\in\T^{N-1}$ for some limit function~$U$. Recall that the functions~$U_{A_n}$ are translated in such a way that
\be\label{normtheta}
\max_{y\in\T^{N-1}}\,U_{A_n}(0,y)=\theta.
\ee
Hence, for all $n\in\N$, $0<U_{A_n}\le\theta$ in $[0,+\infty)\times\T^{N-1}$ and
\be\label{eqlin}
\Delta_yU_{A_n}+A_n^{-2}U_{A_n,xx}+(\gamma_{A_n}-\alpha(y))\,U_{A_n,x}= 0\ \hbox{ in }[0,+\infty)\times\T^{N-1}.
\ee\par
On the other hand, for each $\lambda\in\R$, let $\mu_n(\lambda)$ denote the principal eigenvalue of the operator
$$\Delta_y+A_n^{-2}\lambda^2-\lambda\,(\gamma_{A_n}-\alpha(y))$$
in~$\T^{N-1}$, that is $\mu_n(\lambda)$ is associated with positive eigenfunctions of this operator. Since each function $\mu_n$ is convex and satisfies $\mu_n(0)=0$, $\mu_n'(0)=-\gamma_{A_n}<0$ (because $\alpha$ has zero average over $\T^{N-1}$) and $\mu_n(\lambda)\sim A_n^{-2}\lambda^2$ as $\lambda\to+\infty$, there is a unique $\lambda_n>0$ such that $\mu_n(\lambda_n)=0$. Let $\varphi_n$ be an eigenfunction associated with $\lambda=\lambda_n$, that is
$$\left\{\baa{rcl}
\Delta_y\varphi_n+A_n^{-2}\,\lambda_n^2\,\varphi_n-\lambda_n\,(\gamma_{A_n}-\alpha(y))\,\varphi_n & = & 0\ \hbox{ in }\T^{N-1},\vspace{3pt}\\
\varphi_n & > & 0\ \hbox{ in }\T^{N-1}.\eaa\right.$$
Then $e^{-\lambda_nx}\varphi_n(y)$ solves~(\ref{eqlin}). Up to normalization, one can assume that
\be\label{minphin}
\min_{\T^{N-1}}\varphi_n=\theta
\ee
and it then follows from the maximum principle in half-cylinders as in \cite{v1} (see also \cite{bh1}) that
$$U_{A_n}(x,y)\le e^{-\lambda_nx}\varphi_n(y)\ \hbox{ for all }(x,y)\in[0,+\infty)\times\T^{N-1}.$$
Furthermore, for each $n\in\N$, because of (\ref{infsup}) and $A_n\ge 1$, there holds
$$\mu_n(\lambda)\le\underline{\mu}(\lambda)\ \hbox{ for all }\lambda\ge 0,$$
where $\underline{\mu}(\lambda)$ denotes the principal eigenvalue of the operator $\Delta_y+\lambda^2-\lambda\,(\underline{\gamma}-\alpha(y))$ in~$\T^{N-1}$. As above, there is a unique $\underline{\lambda}>0$ such that $\underline{\mu}(\underline{\lambda})=0$, and $\underline{\mu}(\lambda)<0$ for all $\lambda\in(0,\underline{\lambda})$. As a consequence, $\lambda_n\ge\underline{\lambda}>0$ for all $n\in\N$, whence
\be\label{majoration}
U_{A_n}(x,y)\le e^{-\underline{\lambda}x}\varphi_n(y)\ \hbox{ for all }(x,y)\in[0,+\infty)\times\T^{N-1}.
\ee\par
Let $y_n\in\T^{N-1}$ be such that $\varphi_n(y_n)=\theta$.
 Up to extraction of a subsequence, one can assume that $y_n\to\underline{y}\in\T^{N-1}$ as $n\to+\infty$.  Pick any $\epsilon\in(0,2\theta)$ and define
$$\underline{x}=-\underline{\lambda}^{-1}\ln\Big(\frac{\epsilon}{2\theta}\Big)\ge 0.$$
Because of~(\ref{majoration}), there holds $U_{A_n}(\underline{x},y_n)\le\epsilon/2$ for all $n\in\N$. But since the sequence~$(\|\nabla_yU_{A_n}\|_{L^{\infty}(\R\times\T^{N-1})})_{n\in\N}$ is bounded by Lemma~\ref{lem1}, and since $y_n\to\underline{y}$ as $n\to+\infty$, there are $r>0$ and $n_0\in\N$ such that
$$U_{A_n}(\underline{x},y)\le\epsilon\ \hbox{ for all }|y-\underline{y}|\le r\hbox{ and }n\ge n_0.$$
By monotonicity in $x$, it follows that $U_{A_n}(x,y)\le\epsilon$ for all $x\ge\underline{x}$, $|y-\underline{y}|\le r$ and $n\ge n_0$. By passing to a limit as $n\to +\infty$ as in step 1 (up to extraction of a sequence), the same inequality holds for the limit function~$U$ for almost every $(x,y)$ such that $x\ge\underline{x}$ and $|y-\underline{y}|\le r$. Because of~(\ref{monotone}) and step~2, the constant limit $U(+\infty,\cdot)$ of $U$ as $x\to+\infty$ is then such that
$$0\le U(+\infty,\cdot)\le\epsilon.$$
Since $\epsilon>0$ was arbitrarily small, one concludes that $U(+\infty,\cdot)=0$.\par
Let us now prove that $U(-\infty,\cdot)=1$. First, because of~(\ref{normtheta}), there is a sequence $(y'_n)_{n\in\N}$ in $\T^{N-1}$ such that $U_{A_n}(0,y'_n)=\theta$. Up to extraction of another subsequence, one can assume that $y'_n\to\overline{y}\in\T^{N-1}$ as $n\to+\infty$. As in the previous paragraph, given any $\epsilon\in(0,\theta)$, it follows from Lemma~\ref{lem1} that there are $r'>0$ and $n'_0\in\N$ such that
$$U_{A_n}(0,y)\ge\theta-\epsilon\hbox{ for all }|y-\overline{y}|\le r'\hbox{ and }n\ge n'_0,$$
whence $U_{A_n}(x,y)\ge\theta-\epsilon$ for all $x\le 0$, $|y-\overline{y}|\le r'$ and $n\ge n'_0$. Thus $U(x,y)\ge\theta-\epsilon$ for a.e. such $(x,y)$ and the constant function $U(-\infty,\cdot)$ satisfies $U(-\infty,\cdot)\ge\theta-\epsilon$ for all $\epsilon\in(0,\theta)$. Therefore, $U(-\infty,\cdot)\ge\theta$. Because of~(\ref{f}) and~(\ref{fuinfty}), the constant $U(-\infty,\cdot)$ is either~$\theta$ or~$1$.\par
Assume now that $U(-\infty,\cdot)=\theta$. The monotonicity property~(\ref{monotone}) yields $0\le U\le\theta$ a.e. in $\R\times\T^{N-1}$, whence
\be\label{eqUter}
\Delta_yU+(\gamma-\alpha(y))\,U_x=0\ \hbox{ in }{\mathcal{D}}'(\R\times\T^{N-1}).
\ee
Let $\phi:[0,1]\to\R$ be a non-increasing $C^{\infty}$ function such that $\phi=1$ on $[0,1/3]$ and $\phi=0$ on~$[2/3,1]$. For each $k\in\N$, let $\phi_k$ be the function defined in $\R\times\T^{N-1}$ by
$$\phi_k(x,y)=\left\{\baa{ll}0 & \hbox{if }|x|\ge k+1,\vspace{3pt}\\
1 & \hbox{if }|x|\le k,\vspace{3pt}\\
\phi(|x|-k) & \hbox{if }k<|x|<k+1.\eaa\right.$$
By testing (\ref{eqUter}) against $\phi_k$, one gets that
$$\int_{(0,1)\times\T^{N-1}}\!\!\!\!\!\!\!\!(\gamma-\alpha(y))\,U(x+k,y)\,\phi'(x)\,dx\,dy  -\int_{(0,1)\times\T^{N-1}}\!\!\!\!\!\!\!\!(\gamma-\alpha(y))\,U(x-k-1,y)\,\phi'(1-x)\,dx\,dy=0.$$
Since $U(+\infty,\cdot)=0$, $U(-\infty,\cdot)=\theta$ and the function~$\alpha$ has zero average over $\T^{N-1}$, it follows from Lebesgue's dominated convergence theorem, by passing to the limit as $k\to+\infty$, that~$\theta\,\gamma=0$. But $\theta>0$ and $\gamma>0$ by~(\ref{gamma+}), a contradiction.  So $U(-\infty,\cdot)=1$ and the uniformity of the limit in $y$ (as well as of $U(+\infty,\cdot)=0$) follows from Lemma \ref{lem1} and $U$ being non-increasing in $x$. The proof of Proposition~\ref{proex1} is complete.\hfill$\Box$\break

To complete this section, we prove the existence of solutions~$(\gamma,U)$ of~(\ref{eqU}) when $\theta=0$, for all
$$\gamma\ge\gamma^*_{\infty}=\liminf_{A\to+\infty}\gamma^*_A.$$
Notice also that 
\be\label{2.1a}
0<\gamma^*_{\infty}\le\max_{\T^{N-1}}\alpha
\ee
by~(\ref{limsup}) and~(\ref{infsup}).

\begin{pro}\label{proex2} Assume that $\theta=0$. Then, for all $\gamma\ge\gamma^*_{\infty}$, there exists a solution~$U$ of~$(\ref{eqU})$ such that $\nabla_yU\in L^2(\R\times\T^{N-1})\cap L^{\infty}(\R\times\T^{N-1})$.
\end{pro}

\noindent{\bf{Proof.}} Fix a real number $\gamma\in[\gamma^*_{\infty},+\infty)$. Since $\gamma^*_A$ is the minimal speed for the solutions of~(\ref{eqUA}) for each $A>0$, there exist sequences~$(A_n)_{n\in\N}$ and $(\gamma_n)_{n\in\N}$ of positive real numbers such that
$$A_n\to+\infty\hbox{ and }\gamma_n\to\gamma\hbox{ as }n\to+\infty,$$
and $\gamma_n\ge\gamma^*_{A_n}$ for all $n\in\N$. For each $n\in\N$, let $U_{A_n}$ be a solution of~(\ref{eqUA}) with $A=A_n$ and speed $\gamma_n$ (so hypotheses of Lemma \ref{lem1} are satisfied). Up to a shift in $x$, one can assume
\be\label{norm2}
\int_{(0,1)\times\T^{N-1}}U_{A_n}(x,y)\,dx\,dy=\frac{1}{2}.
\ee
From Lemma~\ref{lem1}, it follows as in step~1 of the proof of Proposition~\ref{proex1} that, up to extraction of a subsequence, the functions~$U_{A_n}$ converge a.e. in $\R\times\T^{N-1}$ to a function~$U\in L^{\infty}(\R\times\T^{N-1})$ such that $0\le U\le 1$ a.e. in $\R\times\T^{N-1}$ and $\nabla_yU\in L^2(\R\times\T^{N-1})\cap L^{\infty}(\R\times\T^{N-1})$. Furthermore,~$U$ fulfills~(\ref{monotone}), and~(\ref{eqUbis}) in $\mathcal{D}'(\R\times\T^{N-1})$. Lastly, from step~2 of the proof of Proposition~\ref{proex1}, the limits $U(\pm\infty,\cdot)$ exist, are constant, and are zeros of the function~$f$ in $[0,1]$. By passing to the limit as $n\to+\infty$ in~(\ref{norm2}), it follows from Lebesgue's dominated convergence theorem that
$$\int_{(0,1)\times\T^{N-1}}U(x,y)\,dx\,dy=\frac{1}{2}.$$
Because of~(\ref{monotone}), one concludes that $0\le U(+\infty,\cdot)\le1/2\le U(-\infty,\cdot)\le 1$.  Since $0$ and $1$ are the only zeros of $f$ in $[0,1]$ we have
$$U(+\infty,\cdot)=0\ \hbox{ and }\ U(-\infty,\cdot)=1,$$
which completes the proof of Proposition~\ref{proex2}.\hfill$\Box$

%%%%%%%%%%%%%%%%%%%%%%%%%%%%%%%%%%%%%%%%
%%%%%%%%%%%%%%%%%%%%%%%%%%%%%%%%%%%%%%%%

\section{Proof of Theorem~\ref{th1}}\label{sec3}

This section is devoted to the completion of the proof of Theorem~\ref{th1}. The existence of a pair~$(\gamma,U)$ sol\-ving~(\ref{eqU}) has already been proven in the previous section. First, in Subsection~\ref{sec31}, we prove its uniqueness, up to shifts in $x$ for~$U$, when  the additional hypothesis~(\ref{hypalpha}) is satisfied. The proof relies on the regularity theory \cite{rs} and a strong maximum principle \cite{sv} for degenerate elliptic and parabolic equations. From this uniqueness property and the Lipschitz continuity of the unique speeds~$c^*(Aq,f)$ with respect to the flow~$Aq$, we deduce in Subsection~\ref{sec32} the existence of~$ \lim_{A\to+\infty} c^*(Aq,f)/A$ with or without~(\ref{hypalpha}).

%%%%%%%%%%%%%%%%%%%%%%%%%%%%%%%%%%%%%%%%

\subsection{Uniqueness of solutions of~(\ref{eqU}) given~(\ref{hypalpha}) and $0<\theta<1$}\label{sec31}

Let~$(\gamma,U)$ and $(\gamma',U')$ be two solutions of \eqref{eqU}.  We will prove that~$\gamma=\gamma'$, and~$U=U'$ after a shift in $x$, under the additional hypothesis \eqref{hypalpha}.

First, because of~(\ref{hypalpha}), it follows from~\cite[Theorem 18(c)]{rs} with $X_j=\partial_{y_j}$ ($j=1,\dots,N-1$) and $X_0=(\gamma-\alpha(y))\partial_x$ (resp.~$X_0=(\gamma'-\alpha(y))\partial_x$) that~$U$ and~$U'$ are actually continuous. In fact, $U$  and $U'$ are slightly smoother than $f$, at least of class $C^{1,\delta}(\R\times\T^{N-1})$), by repeated application of~\cite[Theorem 18(b)]{rs}.

Next,  \cite[Theorem 18(d)]{rs} and the fact that $f(U)\in L^p_{\rm loc}(\R\times\T^{N-1})$  (and so $f(U)$ is in the space $S^p_0(M)$ from \cite{rs}, for any bounded open $M\subseteq \R\times\T^{N-1}$) yield $U\in S^p_2(M)$ for all $p\in(1,\infty)$.  This means that $\Delta_yU, (\gamma-\alpha(y))U_x\in L^p_{\rm loc}(\R\times\T^{N-1})$ for any $p\in(1,\infty)$, so \eqref{eqU} holds in $L^p_{\rm loc}(\R\times\T^{N-1})$ for  $p\in(1,\infty)$.  Thus $LU\le 0$ and $LU'\le 0$ in $L^p_{\rm loc}(\R\times\T^{N-1})$, where
$$L=\Delta_y+(\gamma-\alpha(y))\partial_x.$$  

We will now show that this means that the strong maximum principle from \cite{sv} applies to $-U$ and $-U'$. In order to show they satisfy the hypotheses in \cite{sv}, we will need to approximate them by smooth functions.
Let $\eta:\R\to\R$ be a smooth non-negative function supported in $[-1,1]$, with $\int_{-1}^1\eta(\zeta)d\zeta=1$ and  $\|\eta'\|_\infty\le 10$.  For $\eps>0$ and $\kappa\ge 0$ define the mollifiers 
\[
\eta_{\eps,\kappa}(x,y) :=
\begin{cases}
\eps^{-1}\kappa^{1-N} \eta \left( \displaystyle{\frac x\eps} \right)  \displaystyle{\mathop{\prod}_{j=1}^{N-1}} \eta \left( \displaystyle{\frac {y_j}\kappa} \right) & \kappa>0,
 \\ \eps^{-1} \eta \left( \displaystyle{\frac x\eps} \right) \delta_0(y) & \kappa=0,
\end{cases}
\]
where $\delta_0$ is the delta function at $0\in\R^{N-1}$.  Extend $U$ to a periodic-in-$y$ function on $\R^N$ and let
$$U^{\eps,\kappa} := U*\eta_{\eps,\kappa}.$$
We have $|U^{\eps,\kappa}_x|\le 10\eps^{-1}$ since $|U|\le 1$.  So if $\omega\in C([0,{+\infty}))$ with $\omega(0)=0$ is such that  $|\alpha(y)-\alpha(y')| \le \omega(\kappa)$ whenever $|y-y'|\le \sqrt{N-1}\,\kappa$, then for $\kappa>0$ we have
\[
LU^{\eps,\kappa} = (LU)*\eta_{\eps,\kappa} + \int_{\R^{N-1}} [\alpha(y-\xi)-\alpha(y)]U^{\eps,0}_x(x,y-\xi) \kappa^{1-N}\prod_{j=1}^{N-1} \eta \left( \frac {\xi_j}\kappa \right) \, d\xi \le  10\eps^{-1}\omega(\kappa).
\]
So if we define $\bar U^{\eps,\kappa}  :=U^{\eps,\kappa} - 5\eps^{-1}\omega(\kappa) y_1^2$, then $L\bar U^{\eps,\kappa}\le 0$.  Since also $\bar U^{\eps,\kappa}\in C^2(\R^N)$ for $\kappa>0$, the last claim of the first paragraph of \cite[Section 6]{sv} tells us that $-\bar U^{\eps,\kappa}\in\mathcal H^-_L(\R^N)$, that is, it has the submartingale property relative to $L$ described in that paragraph (with $\bar U^{\eps,\kappa}$ independent of the variable $t$).  Finally, since $U$ is continuous, $\bar U^{\eps,\kappa}\to U^{\eps,0}$ as $\kappa\to 0$ and then $U^{\eps,0}\to U$ as $\eps\to 0$, both locally uniformly in $\R^N$.  This means that $-U^{\eps,0}$ and then $-U$ have the same submartingale property, hence $-U\in\mathcal H^-_L(\R^N)$ (and also $-U\in \mathcal{H}^-_L(\R\times\T^{N-1})$).  Thus the strong maximum principle from \cite[Theorem 6.1]{sv} applies to $-U$ (and to $-U'$), which we will use in the rest of this section.  

We will also need the following.

\begin{lem}\label{lemstrict} If~$(\ref{hypalpha})$ holds and~$(\gamma,U)$ is a solution of \eqref{eqU}, then
\be\label{strict}
\int_{\T^{N-1}} \alpha(y) dy < \gamma<\max_{y\in \T^{N-1}}\,\alpha(y).
\ee
\end{lem}

\noindent{\bf{Proof.}} 
For $a_1<a_2\le b_1<b_2$, let $\phi_{a_1,a_2,b_1,b_2}:\R\to[0,1]$ be a smooth compactly supported test function, equal to 0 on  $(-\infty,a_1]\cup[b_2,+\infty)$, to 1 on $[a_2,b_1]$, increasing on $[a_1,a_2]$, and decreasing on $[b_1,b_2]$.  Test the PDE \eqref{eqU} against $\phi_{a,a+1,b,b+1}$ (extended onto $\R\times\T^{N-1}$ so it is independent of $y$) and take $-a,b\to +\infty$ to obtain
\[
\int_{\T^{N-1}}(\gamma-\alpha(y)) dy = \int_{\R\times\T^{N-1}} f(U(x,y))dxdy>0
\]
(the last inequality is due to $U$ being continuous).  The first inequality in \eqref{strict} now follows.

Next test  against $\phi_{a,a+\epsilon,b,b+1}$, where $a$ is the smallest real number such that
$$\|U\|_{L^\infty([a,+\infty)\times\T^{N-1})}=\theta$$
(recall that $U$ is continuous).  Then take $\epsilon\to 0$ and $b\to{ +\infty}$ to obtain
\[
\int_{\T^{N-1}}(\gamma-\alpha(y)) U(a,y) dy = \int_{[a,+\infty)\times\T^{N-1}} f(U(x,y))dxdy = 0 .
\]
Condition \eqref{hypalpha}, continuity of $U\ge 0$, and $U(a,y_a)=\theta$ for some $y_a$ now  yield the second inequality in \eqref{strict}. 
\hfill$\Box$\break

Now, since $f(0)=f(1)=0$ and 
\be\label{min0max}
\min_{y\in\T^{N-1}}\,(\gamma-\alpha(y))<0<\max_{y\in\T^{N-1}}\,(\gamma-\alpha(y)),
\ee
it follows from the maximum principle for degenerate parabolic equations in~\cite{sv} that $U$ cannot attain the value~$0$ (resp.~$1$) unless being identically equal to~$0$ (resp.~$1$) in~$\R\times\T^{N-1}$.  (Here we apply \cite[Theorem 6.1]{sv} to the functions $V(t,x,y)=V(x,y)=-U(x,y)\le 0$ and $\tilde V(t,x,y)=e^{-\|f'\|_\infty t}(U(x,y)-1)\le 0$, and use that \eqref{min0max} yields $\mathcal{G}(t_0,x_0,y_0)=[t_0,+\infty)\times\R\times\T^{N-1}$ in that theorem.  Note that \cite[Theorem 6.1]{sv} applies to $\tilde V$ just as it does to $V$, because of its continuity and this time the fact that $\tilde V_t+L\tilde V\ge 0$ in $L^p_{\rm loc}(\R^2\times\T^{N-1})$ for $p\in(1,\infty)$.)
But $U(-\infty,\cdot)=1$ and $U(+\infty,\cdot)=0$, so we have
$$0<U<1 \qquad \text{and}\qquad  0<U'<1 \qquad \hbox{ on }\R\times\T^{N-1}.$$\par

Let us assume here that
\be\label{gamma'}
\gamma'\le\gamma,
\ee
which we can ensure by possibly swapping $U$ and $U'$. We shall now slide $U'$ in the $x$-direction and compare it with $U$. We will first prove that $U\le U'(\cdot-h,\cdot)$ in $\R\times\T^{N-1}$ for~$h$ large enough, and then decrease~$h$ up to a critical value $h_*$ so that $U$ and $U'(\cdot-h_*,\cdot)$ will touch.  It will follow from the maximum principle in \cite{sv}  that $\gamma=\gamma'$ and~$U=U'(\cdot-h_*,\cdot)$. 

\begin{lem}\label{lemh0} There exists $h_0\in\R$ such that
$$U(x,y)\le U'(x-h,y)\ \hbox{ for all }h\ge h_0 \hbox{ and all } (x,y)\in\R\times\T^{N-1}.$$
\end{lem}

\noindent{\bf{Proof.}} Extend $f$ by $0$ on~$\R\backslash[0,1]$ (so it is Lipschitz-continuous on~$\R$).  Let $\rho>0$ be such that $f'\le 0$ on $[1-\rho,1]$ and let $M>0$ be such that
\be\label{ineqsUU'}
U\le \theta  \hbox{ in }[M,+\infty)\times\T^{N-1} \qquad\hbox{and}\qquad U'\ge1-\rho \hbox{ in }(-\infty,-M]\times\T^{N-1}.
\ee
We will prove the lemma with $h_0=4M$.

Pick any $h\ge h_0$ and let
\be\label{3.0}
\epsilon_*=\min\Big\{\epsilon\ge 0 \,\Big|\, U(x,y)-\epsilon\le U'(x-h,y)\hbox{ for all }(x,y)\in\R\times\T^{N-1}\Big\}.
\ee
We want to show that $\epsilon_*=0$, so let us assume $\epsilon_*>0$.  Then continuity of $U$ and $U'$, and the limits in \eqref{eqU}, show that the function
$$V(x,y)=U(x,y)-\epsilon_*-U'(x-h,y)\le 0$$
must attain the value 0 at some $(x_0,y_0)$.  Moreover, we have
\be\label{3.1}
\Delta_yV(x,y)+(\gamma-\alpha(y))\,V_x(x,y) =  (\gamma'-\gamma)U'_x(x-h,y) + f(U'(x-h,y))-f(U(x,y)).
\ee

Assume for now that, in the sense of distributions,
\be\label{3.2}
(\gamma'-\gamma) U'_x \ge 0 \qquad\text{on $\R\times\T^{N-1}$.}
\ee
If $x_0\ge 2M$, let $\mathcal G=(M,+\infty)\times\T^{N-1}$.  
Since $U\le\theta$ on $\mathcal G$, we have $LV\ge 0$  on $\mathcal G$ (as a distribution, but then also in $L^p_{\rm loc}(\mathcal G)$ because $LU,LU'\in L^p_{\rm loc}(\R\times\T^{N-1})$).  Then \cite[Theorem~6.1]{sv} (which again applies to $V$ because it is continuous and $LV\ge 0$ in $L^p_{\rm loc}(\mathcal G)$) together with \eqref{min0max} shows $V=0$  on $\mathcal G$, a contradiction with $\epsilon_*>0$ and the limits in \eqref{eqU}.  

If  instead $x_0\le 2M$ and \eqref{3.2} is assumed, let $\mathcal G=(-\infty,3M)\times\T^{N-1}$ and define
\be\label{3.3}
\tilde V(t,x,y)=e^{-\|f'\|_\infty t} V(x,y)\le 0.
\ee
Since $U'(\cdot-h,\cdot)\ge 1- \rho$ on $\mathcal G$, we have 
\[
f(U'(x-h,y))-f(U(x,y))\ge f(U'(x-h,y)+\epsilon_*)-f(U(x,y))\ge -\|f'\|_\infty |V(x,y)|
\]
there. Thus
\be\label{3.4}
\tilde V_t + L\tilde V \ge 0 
\ee
in $L^p_{\rm loc}(\R\times\mathcal G)$ and \cite[Theorem~6.1]{sv} again shows $\tilde V=0$  on $\R\times\mathcal G$, a contradiction.  

It follows that $\epsilon_*=0$, so the lemma is proved under the assumption \eqref{3.2} (e.g., if $U=U'$, since then $\gamma=\gamma'$).  The validity of \eqref{3.2} will be obtained at the end of the proof of the next lemma.  \hfill$\Box$\break

Next let
\be\label{defh*}
h_*=\min\Big\{h'\in\R \,\Big|\, U(x,y)\le U'(x-h,y)\hbox{ for all }(x,y)\in\R\times\T^{N-1}\hbox{ and }h\ge h'\Big\}.
\ee
If we assume \eqref{3.2} (which we shall prove to be true always), then Lemma~\ref{lemh0} {and the limits in~\eqref{eqU}} show that $h_*$ is well defined. 

\begin{lem}\label{lemgamma} We have $\gamma=\gamma'$ and 
$$U(x,y)=U'(x-h_*,y)\hbox{ for all }(x,y)\in\R\times\T^{N-1}.$$
\end{lem}

\noindent{\bf{Proof.}} Assume again \eqref{3.2}, so $h_*$ is well defined, and let
$$V(x,y)=U(x,y)-U'(x-h_*,y)\le 0.$$
Then \eqref{3.1} and \eqref{3.2} show that $\tilde V$ from \eqref{3.3}  satisfies \eqref{3.4}  on $\R^2\times\T^{N-1}$.  If $V(x_0,y_0)=0$ for some $(x_0,y_0)$, then \cite[Theorem~6.1]{sv} again implies $\tilde V=0$ on $\R^2\times\T^{N-1}$, as was to be proved.  

Let us therefore assume $V<0$ on $\R\times\T^{N-1}$ and pick $\eta\in(0,h_*]$ such that
$$U(x,y)< U'(x-h,y)\hbox{ for all }h\in[h_*-\eta,h_*] \text{ and all }(x,y)\in[-2M,2M]\times\T^{N-1}.$$
This can be done because $U$ and $U'$ are continuous.  For any such $h$ define $\epsilon_*$ as in \eqref{3.0}.  It follows that $U(x_0,y_0)-\epsilon_*-U'(x_0-h,y_0)=0$ for some $(x_0,y_0)$ with $|x_0|>2M$.  As in the proof of the last lemma, this, $U\le\theta$ on $(M,+\infty)\times\T^{N-1}$, and $U'(\cdot-h,\cdot)\ge 1-\rho$ on $(-\infty,-M)\times\T^{N-1}$ imply $U(x,y)-\epsilon_*-U'(x-h,y)=0$ on one of these two sets. In either case we must have $\epsilon_*=0$. This means that~$h_*$ in~(\ref{defh*}) can be decreased by $\eta>0$, a contradiction with the assumption that it is minimal.

Thus we have proved the lemma assuming \eqref{3.2}.  But \eqref{3.2} is satisfied if $U=U'$, so the definition of $h_*$,  the limits in~\eqref{eqU}, and the lemma yield $h^*=0$ in this case and $U'_x\le 0$. Therefore \eqref{3.2} holds a priori (after assuming \eqref{gamma'}), so this and the previous lemma also hold as stated. \hfill$\Box$\break

Finally the function $U_x\le 0$ (which is of class $C^{\delta}(\R\times\T^{N-1})$, as we have showed earlier)  is by~\eqref{eqU} a distributional solution of
$$LU_x=-f'(U(x,y))\,U_x.$$
As at the beginning of this subsection, it follows from~\cite[Theorem 18(d)]{rs} that $\Delta_yU_x,(\gamma-\alpha(y))U_{xx}\in L^p_{\rm loc}(\R\times\T^{N-1})$ for any $p\in(1,\infty)$, so the above equation holds in $L^p_{\rm loc}(\R\times\T^{N-1})$ for any $p\in(1,\infty)$. Since $U_x$ cannot be identically equal to $0$ in $\R\times\T^{N-1}$,~\eqref{min0max} and~\cite[Theorem 6.1]{sv} imply that $U_x<0$ everywhere in $\R\times\T^{N-1}$.

%%%%%%%%%%%%%%%%%%%%%%%%%%%%%%%%%%%%%%%%

\subsection{Existence of the limit~(\ref{gamma}) when $0<\theta<1$}\label{sec32}

The only part of Theorem~\ref{th1} left to prove is the existence of the limit~(\ref{gamma}). If (\ref{hypalpha}) holds, then it follows from the proof of Proposition~\ref{proex1} and the uniqueness result in Lemma~\ref{lemgamma} that the bounded family $(\gamma^*_A)_{A\ge 1}$ has only one limiting value as $A\to+\infty$. In other words, the limit
\be\label{3.5}
\gamma^*(q,f):=\lim_{A\to+\infty}\gamma^*_A=\lim_{A\to+\infty}\frac{c^*(Aq,f)}{A}
\ee
exists. Furthermore, $\gamma^*(q,f)$ satisfies~(\ref{strict}) in this case.

Now consider general $\alpha$, pick  any $\epsilon>0$, let~$\overline{\alpha}$ satisfying \eqref{hypalpha} be such that $\alpha\le \overline{\alpha}\le \alpha+\epsilon$, and let
$$\overline{q}(x,y)=(\overline{\alpha}(y),0,\cdots,0).$$
Then for any $A\in\R$ we have
$$ c^*(A\overline{q},f)-A\epsilon\le  c^*(Aq,f)\le c^*(A\overline{q},f)$$
because the speed of traveling fronts for problem~(\ref{equ}) is $1$-Lipschitz with respect to the~$L^{\infty}$~norm of the flow, see~\cite{bhbook,bn2}.  Since \eqref{3.5} shows that $c^*(A\overline{q},f)/A$ converges as $A\to +\infty$, and since $\epsilon>0$ was arbitrary, existence of the limit \eqref{gamma} follows.  Moreover \eqref{strict} holds for $\overline\alpha$ and $\gamma^*(\overline q,f)$, so we obtain
\[
\int_{\T^{N-1}} \alpha(y) dy \le \gamma^*(q,f) \le \max_{y\in \T^{N-1}}\,\alpha(y),
\]
with the first inequality strict unless $\alpha$ is constant by \cite{kr}.  The proof of Theorem \ref{th1} is finished.\hfill$\Box$\break

\begin{rem}{\rm When the flow~$\alpha$ satisfies the non-degeneracy condition~$(\ref{hypalpha})$, then the solutions~$(\gamma,U)$ of~$(\ref{eqU})$ are unique, up to shifts in~$x$ for~$U$. Without assumption~$(\ref{hypalpha})$, this uniqueness property is not clear in the general case, for which the solution~$U$ obtained in Proposition~$\ref{proex1}$ is just a distributional solution which may not be of class $C^{1,\delta}(\R\times\T^{N-1})$. Furthermore, the strict inequality $\gamma^*(q,f)<\max_{\T^{N-1}}\alpha$, which holds under assumption~$(\ref{hypalpha})$, may not be true in the general case, when $\alpha$ is flat and equal to its maximum on a sufficiently large region.}
\end{rem}

\begin{rem}{\rm Since the speed of traveling fronts for problem~$(\ref{equ})$ is $1$-Lipschitz with respect to the~$L^{\infty}$~norm of the flow, it follows from~$(\ref{gamma})$ that, if $q_{\alpha}(x,y)=(\alpha(y),0,\ldots,0)$ and $q_{\beta}(x,y)=(\beta(y),0,\ldots,0)$, then
$$|\gamma^*(q_{\alpha},f)-\gamma^*(q_{\beta},f)|\le\|\alpha-\beta\|_{L^{\infty}(\T^{N-1})}.$$}
\end{rem}

%%%%%%%%%%%%%%%%%%%%%%%%%%%%%%%%%%%%%%%%
%%%%%%%%%%%%%%%%%%%%%%%%%%%%%%%%%%%%%%%%

\section{Proof of Theorem~\ref{th2}}\label{sec4}

The existence result was proved in Section \ref{sec2} with
$$\gamma^*(q,f):=\gamma^*_\infty = \liminf_{A\to\infty} \gamma^*_{A}.$$
It remains to prove \eqref{gamma}  and that $\gamma^*(q,f)$ equals the limit in the statement of the theorem (for general $\alpha$) and the claims in the last sentence (for $\alpha$ satisfying \eqref{hypalpha}). In fact, it will be sufficient to assume \eqref{hypalpha} in the rest of the proof because \eqref{gamma} as well as the $\theta'\to 0$ limit in the general case then follow as in Subsection \ref{sec32}. 

So let us assume \eqref{hypalpha}. Recall that $f_{\theta'}(u):= f(u)\chi(u/{\theta'})$ for  ${\theta'}\in(0,1/4]$ and let $\widetilde{f}_{\theta'}$ be the function which equals $f_{\theta'}$ on $[0,1]$ and equals $0$ on $[-\theta',0)$. For $A\ge 1$ and $\theta'\in(0,1/4]$, let $\widetilde{\gamma}^*_{A,\theta'}$ (resp. $\gamma^*_{A,\theta'}$) be the unique speed for which there is a (unique up to shifts) solution $\widetilde{U}^{A,\theta'}$ (resp. $U^{A,\theta'}$) of~\eqref{eqUA} with $\widetilde{f}_{\theta'}$ (resp. $f_{\theta'}$) in place of $f$, but with $\widetilde{U}^{A,\theta'}$ instead satisfying $-\theta'<\widetilde{U}^{A,\theta'}<1$ and $\widetilde{U}^{A,\theta'}(+\infty,\cdot)\equiv-\theta'$.  As in Theorem~\ref{th1}, $\gamma^*(q,\widetilde{f}_{\theta'})=\lim_{A\to\infty}\widetilde{\gamma}^*_{A,\theta'}$ exists and is the unique speed for which there is a (unique up to shifts) solution $\widetilde{U}^{\theta'}$ of~\eqref{eqU} with $\widetilde{f}_{\theta'}$ in place of $f$, and with $-\theta'<\widetilde{U}^{\theta'}<1$ and $\widetilde{U}^{\theta'}(+\infty,\cdot)\equiv-\theta'$.  

It follows from~\cite{bh1} that $\widetilde{\gamma}^*_{A,\theta'}<\gamma^*_{A,\theta'}<\gamma^*_A$ for all $A\ge 1$ and $\theta'\in(0,1/4]$, so
\[
\gamma^*(q,\widetilde{f}_{\theta'})\le\gamma^*(q,f_{\theta'})\le \gamma^*(q,f) .
\]
 Since the maps $\theta'\mapsto\widetilde{\gamma}^*_{A,\theta'}$ and $\theta'\mapsto\gamma^*_{A,\theta'}$ are decreasing on $(0,1/4]$ by~\cite{bh1}, the limits $\lim_{\theta'\to0}\gamma^*(q,\widetilde{f}_{\theta'})$ and $\lim_{\theta'\to0}\gamma^*(q,f_{\theta'})$ exist in $\R$, and
\be\label{4.2bis}
\widetilde{\gamma}:=\lim_{\theta'\to0}\gamma^*(q,\widetilde{f}_{\theta'})\le\gamma:=\lim_{\theta'\to0}\gamma^*(q,f_{\theta'})\le \gamma^*(q,f) \quad (=\liminf_{A\to\infty}\gamma^*_A).
\ee
\noindent{}To obtain \eqref{gamma} and $\gamma^*(q,f) = \lim_{A\to\infty} \gamma^*_{A} = \lim_{{\theta'}\to 0} \gamma^*(q,f_{\theta'})$, we thus only need to prove
\be\label{4.3}
\limsup_{A\to\infty} \gamma^*_{A} \le\widetilde{\gamma}.
\ee

We first note that since front speed is monotone in $f$ and any $f$ is dominated by some KPP $g$ (i.e., with $g(u)\le g'(0)u$), (\ref{hypalpha}) and~(\ref{gamma*kpp}) immediately imply
\be\label{4.4}
\limsup_{A\to\infty} \gamma^*_{A} < \max_{\T^{n-1}} \alpha.
\ee
On the other hand, since $\alpha$ is not constant, Theorem~\ref{th1} yields $\gamma^*(q,\widetilde{f}_{\theta'})>\int_{\T^{N-1}}\alpha(y)dy$ for any $\theta'\in(0,1/4]$, whence
\be\label{4.4bis}
\int_{\T^{N-1}}\alpha(y)\,dy<\widetilde{\gamma}.
\ee

Let us now prove~\eqref{4.3}. With the above notations, Lemma~\ref{lem1} and the boundedness of $\widetilde{\gamma}^*_{A,\theta'}$ (because $\int_{\T^{N-1}}\alpha(y)dy<\widetilde{\gamma}^*_{A,\theta'}<\max_{\T^{N-1}}\alpha+A^{-1}c^*(0,f)$ by~\cite{bh1}) imply that the quantities~$\|\widetilde{U}^{A,\theta'}_x\|_{L^1(\R\times\T^{N-1})}\,(=1)$ and~$\|\nabla_y\widetilde{U}^{A,\theta'}\|_{L^2(\R\times\T^{N-1})}$ are bounded independently of $A\ge 1$ and $\theta'\in(0,1/4]$. Furthermore, since $\sup_{\theta'\in(0,1/4]}\|\widetilde{f}_{\theta'}\|_{C^1([-\theta',1])}<+\infty$, it follows from~\cite{bh2} that the quantities~$\|\nabla_y\widetilde{U}^{A,\theta'}\|_{L^{\infty}(\R\times\T^{N-1})}$ are bounded independently of $A\ge 1$ and $\theta'\in(0,1/4]$. Therefore, from the characterization of~$\widetilde{U}^{\theta'}$ in Theorem~\ref{th1}, there holds
\be\label{suptheta'}
\sup_{\theta'\in(0,1/4]}\Big(\|\widetilde{U}^{\theta'}_x\|_{L^1(\R\times\T^{N-1})}+\|\nabla_y\widetilde{U}^{\theta'}\|_{L^2(\R\times\T^{N-1})}+\|\nabla_y\widetilde{U}^{\theta'}\|_{L^{\infty}(\R\times\T^{N-1})}\Big)<+\infty.
\ee
Up to shifts, one can assume without loss of generality that $\max_{\T^{N-1}}\widetilde{U}^{\theta'}(0,\cdot)=1/2$ and by Sobolev embedding one can choose a sequence $\theta_n\to0$ such that $\widetilde{U}^{\theta_n}(x,y)\to\widetilde{U}(x,y)\in[0,1]$ for almost every $(x,y)\in\R\times\T^{N-1}$. The function $\widetilde{U}$ is a distributional solution of~(\ref{eqU}), with $\widetilde{\gamma}$ in place of $\gamma$, such that $\widetilde{U}(\cdot+h,\cdot)\le\widetilde{U}$ a.e. in $\R\times\T^{N-1}$ for all $h\ge 0$ and $\nabla_y\widetilde{U}\in L^2(\R\times\T^{N-1})\cap L^{\infty}(\R\times\T^{N-1})$. Furthermore, (\ref{suptheta'}), $\widetilde{U}^{\theta'}_x\le 0$, and $\max_{\T^{N-1}}\widetilde{U}^{\theta'}(0,\cdot)=1/2$ imply that $\pm\|\widetilde{U}\|_{L^{\infty}((\pm\eps,+\infty)\times\T^{N-1})}\le \pm 1/2$ for any $\eps>0$. The $x\to\pm\infty$ limits $\widetilde{U}^{\pm}(y)$ of $\widetilde{U}$ are uniform since $\nabla_y\widetilde{U}\in L^{\infty}(\R\times\T^{N-1})$ and they must be $0$ and $1$ because they solve $\Delta_y\widetilde{U}^{\pm}+f(\widetilde{U}^{\pm})=0$ with $f>0$ on $(0,1)$.

Since $f\in C^{1,1}$, it follows from~\cite[Theorem 18(c,b)]{rs} that $\widetilde U$ is actually a classical solution of class~$C^2(\R\times\T^{N-1})$, with uniformly bounded first and second derivatives (this is true for any solution, proving the first claim in the last sentence of the Theorem).  We now let $V:=-\widetilde U_x\ge 0$ and obtain from \eqref{eqU} for $\widetilde U$,
\be\label{4.5}
\Delta_y V+(\widetilde{\gamma} - \alpha(y))\,V_x=-f'(\widetilde U(x,y))V.
\ee
Observe that \eqref{4.2bis},~\eqref{4.4}, and~\eqref{4.4bis} yield
\be\label{min0maxbis}
\min_{y\in\T^{N-1}}\big(\widetilde{\gamma}-\alpha(y)\big)<0<\max_{y\in\T^{N-1}}\big(\widetilde{\gamma}-\alpha(y)\big).
\ee
Since $f'(\widetilde U(x,y))$ is bounded and Lipschitz, we obtain from \cite[Theorem 18(c)]{rs} and subsequent repeated application of \cite[Theorem 18(b)]{rs} that for each $s>0$ there is $ C_s>0$ such that
\[
|V_x(x,y)| \le C_s \|V\|_{L^\infty(B_s(x,y))}\ \hbox{ for all }(x,y)\in\R\times\T^{N-1}
\]
(here $B_s(x,y)$ is the ball of radius $s$ centered at $(x,y)$).  Since $V\ge 0$ satisfies \eqref{4.5} and \eqref{min0maxbis}, the Harnack inequality from \cite{hz2} holds for $V$, so there is $C_s'>0$ such that
\[
 \|V\|_{L^\infty([x-s,x+s]\times\T^{N-1})} \le C_s' V(x,y)\ \hbox{ for all }(x,y)\in\R\times\T^{N-1}.
\]
Thus for any fixed $s>0$ and $C:=C_sC_s'>0$,
\[
|V_x(x,y)| \le C V(x,y)\ \hbox{ for all }(x,y)\in\R\times\T^{N-1}.
\]
But this means that
\[
\Delta_y \widetilde U+ A^{-2} \widetilde U_{xx} + ( \widetilde{\gamma}+CA^{-2}-\alpha(y))\,\widetilde U_x+f(\widetilde U)\le 0.
\]
On the other hand, as for the function $U$ in Theorem~\ref{th1} under assumption~\eqref{hypalpha},~\cite[Theorem 6.1]{sv} and~\eqref{min0maxbis} imply that $\widetilde{U}_x<0$ everywhere in $\R\times\T^{N-1}$. Since
$$\gamma^*_A=\min_{w\in C^2(\R\times\T^{N-1}),\,w_x<0,\,w(-\infty,\cdot)=1,\,w(+\infty,\cdot)=0}\,\sup_{(x,y)\in\R\times\T^{N-1}}\left(\frac{\Delta_yw+A^{-2}w_{xx}+f(w)}{-w_x}+\alpha(y)\right),$$
as follows from the same arguments as the ones used in~\cite{h1} in the case of infinite cylinders with bounded cross sections and Neumann boundary conditions, this implies $\gamma^*_{A}\le \widetilde{\gamma}+CA^{-2}$.  It follows that \eqref{4.3} holds, and thus also \eqref{gamma} and
$$\gamma^*(q,f)=\lim_{\theta'\to0}\gamma^*(q,f_{\theta'})=\lim_{\theta'\to0}\gamma^*(q,\widetilde{f}_{\theta'}).$$

Finally, assume that $U'$ is any solution of \eqref{eqU} with $\gamma'<\gamma^*(q,f)$ and let $U:=\widetilde{U}^{\theta_n}$ for $\theta_n$ as above and small enough so that $\gamma:=\gamma^*(q,\widetilde{f}_{\theta_n})>\gamma'$.  Then the argument after (3.25) in Subsection 3.1 applied to $U,U'$ and with $\theta$ replaced by $0$ in~\eqref{ineqsUU'}  yields a contradiction in the same way as before, with \eqref{3.1} and  \eqref{3.2} 
replaced by 
\[
\Delta_yV(x,y)+(\gamma'-\alpha(y))\,V_x(x,y)  =  (\gamma'-\gamma)U_x(x-h,y) + f(U'(x-h,y))-f(U(x,y))
\]
and by $(\gamma'-\gamma) U_x \ge 0$, and considering the operator $L'$ from the right-hand side of the above equation instead of $L$.
(Note that the upper bound in \eqref{strict} holds for $\gamma'$ because it holds for $\gamma$ and the lower bound in~\eqref{strict} follows from the first part of the proof of Lemma~\ref{lemstrict} which only uses that $f\ge\not\equiv 0$ on $[0,1]$, so \eqref{min0max} holds for $\gamma'$. Moreover,
\[
\lim_{x\to+\infty} U(x,y)\equiv -\theta_n < 0 \equiv \lim_{x\to+\infty} U'(x,y) 
\]
actually makes the proof even easier.)   Thus we must have $\gamma'\ge\gamma^*(q,f)$ and the proof of Theorem \ref{th2} is finished.\hfill$\Box$

%%%%%%%%%%%%%%%%%%%%%%%%%%%%%%%%%%%%%%%%
%%%%%%%%%%%%%%%%%%%%%%%%%%%%%%%%%%%%%%%%

\end{document}